\documentclass[11pt, reqno]{amsart}
\usepackage{enumitem}
\usepackage{amsthm,amsmath,amssymb}

\usepackage{mathrsfs}

\numberwithin{equation}{section}

\newtheorem{theorem}{Theorem}[section]
\newtheorem{definition}[theorem]{Definition}
\theoremstyle{plain}
\newtheorem{proposition}[theorem]{Proposition}
\newtheorem{lemma}[theorem]{Lemma}
\newtheorem{corollary}[theorem]{Corollary}
\newtheorem{remark}[theorem]{Remark}

\def\al{\aligned}
\def\eal{\endaligned}
\def\be{\begin{equation}}
	\def\ee{\end{equation}}

\def\al{\aligned}

\def\d{\nabla}

\DeclareMathOperator{\vol}{Vol}

\numberwithin{equation}{section}
\begin{document}
\title[]{Eigenvalue estimates for Beltrami-Laplacian under Bakry-\'Emery Ricci curvature condition}
\author{Ling Wu, XingYu Song  and Meng Zhu}
\address{School of Mathematical Sciences and Shanghai Key Laboratory of PMMP, East China Normal University, Shanghai 200241, China}
\email{ 1059662403@qq.com, 1137455189@qq.com, mzhu@math.ecnu.edu.cn}
\date{}

\begin{abstract}
On closed Riemannian manifolds with Bakry-\'Emery Ricci curvature bounded from below and bounded gradient of the potential function, we obtain lower bounds for all positive eigenvalues of the Beltrami-Laplacian instead of the drifted Laplacian. The lower bound of the $k$th eigenvalue depends on $k$, Bakry-\'Emery Ricci curvature lower bound, the gradient bound of the potential function, and the dimension and diameter upper bound of the manifold, but the volume of the manifold is not involved. Especially, these results apply to closed manifolds with Ricci curvature bounded from below.
	
\end{abstract}
\maketitle

\section{Introduction}

Let $(M, g)$ be a Riemannian manifold, $f$ a smooth function on $M$. The Bakry-\'Emery Ricci curvature tensor $Ric+Hess\, f$, first introduced in \cite{BE}, is a natural generalization of the classical Ricci curvature tensor (the case where $f$ is a constant). Here, $Ric$ and $Hess\, f$ represent the Ricci curvature tensor and the hessian of $f$, respectively.

Bakry-\'Emery Ricci curvature being bounded below is the concept of ``Ricci curvature bounded below" for smooth metric space $(M, g, e^{-f}dV)$, namely, $M$ equipped with the distance induced by $g$ and measure $e^{-f}dV$, where $dV$ is the volume element. It can also be extended to general metric measure spaces and used to study Ricci limit spaces (see e.g. \cite{St1}, \cite{St2}, \cite{LV}). Moreover, manifolds with constant Bakry-\'Emery Ricci curvature are so called Ricci solitons, which play a crucial role in the singularity analysis of the Ricci flow (see e.g. \cite{Per}, \cite{TZ}, \cite{CW}, \cite{Bam}). Therefore, the question that whether the results for manifolds with Ricci curvature bounded below can also be established when Bakry-\'Emery Ricci curvature is bounded below has drawn a lot of attention.

In this paper, we study the eigenvalue estimates of Beltrami Laplacian $\Delta$ on closed manifolds. The basic assumptions are that $(M^m, g)$ is an $m$-dimensional closed Riemannian manifold with
\be\label{basic assumption1}
Ric+Hess\,f\geq -K g,
\ee
and
\be\label{basic assumption2}
|\d f|\leq L,
\ee
where $\d f$ is the gradient of $f$, and $K$ and $L$ are nonnegative constants.

On manifolds with Ricci curvature bounded below, there have been numerous results on eigenvalue estimates (see e.g. \cite{Lich}, \cite{Ob}, \cite{Cheeger}, \cite{Cheng}, \cite{LY},  \cite{ZY}, \cite{Li}). For manifolds with Bakry-\'Emery Ricci curvature bounded from below, normally the weighted measure $e^{-f}dV$ is considered, and the corresponding self-adjoint Laplace operator is the drifted Laplacian $\Delta_f=\Delta- \d f\cdot \d$. Under the assumptions \eqref{basic assumption1} and \eqref{basic assumption2}, Munteanu-Wang \cite{MW}, Su-Zhang \cite{SZ}, and  Wu \cite{Wu} independently obtained a Cheng type upper bound for the first positive eigenvalue of $\Delta_f$. On the other hand, Charalambous-Lu-Rowlett \cite{CLR} proved lower bound estimates for all positive eigenvalues of $\Delta_f$. An eigenvalue comparison for the first positive eigenvalue of $\Delta_f$ is also given in \cite{BQ} and \cite{AN}.

Different from the above setting, we consider here the standard measure $dV$ and Beltrami-Laplacian $\Delta$ under conditions \eqref{basic assumption1} and \eqref{basic assumption2}. A main difficulty rising in this case is that the hessian of $f$ does not appear in the Bochner formula for $\Delta$, as opposed to the Bochner formula for $\Delta_f$. Thus, to utilize the lower boundedness of the Bakry-\'Emery Ricci curvature, we need to manually add $Hess\,f$, which causes an extra bad term $-Hess\,f(\d \cdot, \d \cdot)$. By using integration by parts and Moser iteration, we are able to overcome this difficulty.

Denote the eigenvalues of $\Delta$ by $0=\lambda_0<\lambda_1\leq \lambda_2\leq \cdots\leq \lambda_k\leq \cdots$, we derive lower bounds for all $\lambda_k$'s. More precisely, we show that

 \begin{theorem}\label{main theorem}
 Let $(M^m,g)$ be an $m$-dimensional closed Riemannian manifold. Assume that conditions \eqref{basic assumption1} and \eqref{basic assumption2} are satisfied. Then\\
(1) we have
\begin{equation}\label{lower bound of 1}
 	\lambda_1\ge c_0;
\end{equation}
(2) for $m \ge 3$,
\begin{equation}\label{lower bound of k}
	\lambda_k\ge c_1k^{\frac{2}{m}}, \ \forall k\geq 2,
\end{equation}
and for $m=2$,
\begin{equation}\label{lower bound of k m=2}
	\lambda_k\ge c_2k^{\frac{1}{2}}, \ \forall k\geq 2.
\end{equation}
Here $c_0$, $c_1$ and $c_2$ are constants depending on $m$, $K,\ L$, and the upper bound $D$ of the diameter of $M$.

\end{theorem}

We prove (1) and (2) of Theorem \ref{main theorem} separately in sections 2 and 3 (see Theorem \ref{thm 1} and Theorem \ref{thm n}), where explicit expressions of $c_0$ , $c_1$ and $c_2$ can also be found. In section 2, we establish the estimate \eqref{lower bound of 1} by finding a lower bound of Cheeger's isoperimetric constant $IN_1(M)$. Actually, we obtain lower bound for the general isoperimetric constant $IN_{\alpha}(M)$, $\alpha>0$, defined in \cite{Li}. The proof follows a method of Dai-Wei-Zhang \cite{DWZ} and uses the volume comparison result of Q. Zhang and the third author \cite{ZZ}. In section 3, following the method in \cite{WZ} (see also \cite{LZZ}), estimates \eqref{lower bound of k} and \eqref{lower bound of k m=2} are proved by using \eqref{lower bound of 1} and gradient estimates for eigenfunctions. The gradient estimates are done by Moser iteration, in which the Sobolev inequality required comes from the isoperimetric constant estimate in section 2.\\


\section{Isoperimetric constant estimate and lower bound of $\lambda_1$ }
In this section, we prove part (1) of Theorem \ref{main theorem}. According to \cite{Cheeger}, it suffices to bound Cheeger's isoperimetric constant from below. Firstly, let us recall the definitions of isoperimetric constants. We adapt the notations and definitions in \cite{Li}.
\begin{definition}
Let $(M,g)$ be a compact Riemannian manifold (with or without boundary). For $\alpha>0$, The Neumann $\alpha$-isoperimetric constant of M is defined by
$$	IN_\alpha(M)=\inf_{\substack{\partial\Omega_1=H=\partial\Omega_2 \\ M=\Omega_1\cup H \cup \Omega_2}}\frac{\vol(H)}{\min\{{\vol(\Omega_1),\vol(\Omega_2)}\}^{\frac{1}{\alpha}}},$$
\\
where the infimum is taken over all hypersurfaces $H$ dividing  $M$ into two parts, denoted by $\Omega_1$ and $\Omega_2$, and  $\vol(\cdot)$ denotes the volume of a region.
\end{definition}

In \cite{Cheeger}, Cheeger showed that
\begin{lemma}\label{Cheeger}
Let $(M,g)$ be a closed Riemannian manifold. Then
 \[\lambda_1 \ge \frac{IN_1(M)^2}{4}.\]
\end{lemma}

Thus, one can get a lower bound of $\lambda_1$ by bounding $IN_1(M)$ from below. As indicated in \cite{DWZ}, this can be done by using the method therein. For completeness, we state the result and also include the proof in the following.


\begin{theorem}\label{isoperimetric estimate}
Let $(M^m,g)$ be an $m$-dimensional complete Riemannian manifold, $m\geq 2$. Assume that \eqref{basic assumption1} and \eqref{basic assumption2} are satisfied. Let $\Omega$ be a bounded convex domain in $M$. Then
for $1\le\alpha \le \frac{m}{m-1},$ we have
	\begin{equation}
	IN_\alpha(\Omega)\ge d^{-1}2^{-2m-1}5^{-m}e^{-(24-\frac{2}{\alpha})Ld-(104-\frac{1}{\alpha})Kd^2} \vol(\Omega)^{1-\frac{1}{\alpha}},
		\end{equation}
	and for $0<\alpha<1,$ we have
	\begin{equation}
	IN_{\alpha}(\Omega)\ge  d^{-1}2^{-2m-1}5^{-m}e^{-22Ld-103Kd^2}\vol(\Omega)^{1-\frac{1}{\alpha}},
		\end{equation}
	where $d$ is the diameter of the domain $\Omega$.

In particular, if $M$ is closed, then
	\begin{equation}
	IN_1(M)\ge D^{-1}2^{-2m-1}5^{-m}e^{-22LD-103K D^2},
		\end{equation}
	and
	\begin{equation}
	IN_{\frac{m}{m-1}}(M)\ge D^{-1}2^{-2m-1}5^{-m}e^{-(22+\frac{2}{m})LD-(103+\frac{1}{m})K D^2}\vol(M)^{\frac{1}{m}},
		\end{equation}	
where $D$ is an upper bound of the diameter of $M$.
	
\end{theorem}

Before starting the proof of Theorem \ref{isoperimetric estimate}, let us present some results needed. First of all, Q. Zhang and the third author \cite{ZZ} proved a volume comparison theorem for manifolds satisfying \eqref{basic assumption1} and \eqref{basic assumption2}.

\begin{theorem}[\cite{ZZ}]\label{volume element comparison}
	Let $(M^m, g)$ be an $m$-dimensional complete Riemannian manifold.  Suppose that $Ric+\frac{1}{2}\mathscr{L}_Vg \ge -Kg$ for some constant $K\ge0$ and smooth vector field $V$ with $|V|\le L$, where $\mathscr{L}_V$ means the Lie derivative in the direction of $V$. Then the following conclusions are true.
	\\
	(a)Let $A(s,\theta)$ denote the volume element of the metric $g$ on M in geodesic polar coordinates. Then for any $0< s_1 <s_2$, we have
\begin{equation}\label{AC}
	\frac{A(s_2,\theta)}{s_2^{m-1}}\le e^{2Ls_2+Ks_2^2} \frac{A(s_1,\theta)}{s_1^{m-1}}.
\end{equation}
\\
(b)For any $0<r_1<r_2$, we have
	\begin{equation} \label{VC}
	\frac{\vol(B_{r_2}(x))}{r_2^m}\le e^{[K(r_2^2-r_1^2)+2L(r_2-r_1)]}\frac{\vol(B_{r_1}(x))}{r_1^m},
\end{equation}
where $B_r(x)$ is the geodesic ball centered at $x\in M$ with radius $r$.
\end{theorem}

\begin{remark}
When $V=\d f$, the assumptions in the above Theorem become \eqref{basic assumption1} and \eqref{basic assumption2}.
\end{remark}

Next, we need the following lemma by Gromov.

\begin{lemma}[\cite{Gro}]\label{Gromov}
	Let $(M^m,g)$ be a complete Riemannian manifold. Let $\Omega$ be a convex domain in $M$, and $H$ a hypersurface dividing $\Omega$ into two parts $\Omega_1,\Omega_2$. For any Borel subsets $W_i \subset \Omega_i,i=1,2$, there exists an $x_1$
	in one of $W_i$, say $W_1$, and a subset $W$ in the other part $W_2$, such that
	\begin{equation}
	\vol(W) \ge \frac{1}{2}\vol(W_2),
	\end{equation}
	and for any $x_2\in W$, there is a unique minimal geodesic $\gamma_{x_1, x_2}$ between  $x_1$ and $x_2$ which intersects $H$ at some $z$ with
		\begin{equation}
	dist(x_1,z)\ge dist(x_2,z),
		\end{equation}
	where $dist(x_1,z)$ denotes the distance between $x_1$ and $z$.	
\end{lemma}

Combining Theorem \ref{volume element comparison} and Lemma \ref{Gromov}, we get

\begin{lemma}
Let $H,W$ and $x_1$ be as in Lemma \ref{Gromov}. Then
	\begin{equation}
 \vol(W)\le D_12^{m-1}e^{4LD_1+4KD_1^2}\vol(H^{'}),
 	\end{equation}
where $D_1=\sup_{x\in W} dist(x_1,x)$, and $H^{'}$ is the set of intersection points with $H$ of geodesics $\gamma_{x_1,x} $ for all $x \in W$.
\end{lemma}

\proof Let $S_{x_1}$ be the set of unit tangent vectors of $M$ at $x_1$, and $\Gamma \subset S_{x_1} $ the subset of vectors $\theta$ such that $\gamma_{\theta} =\gamma_{x_1,x_2}$ for some $x_2\in W$. The volume element of the metric $g$ is written as $dV=A(\theta,t)d\theta \wedge dt$ in polar coordinates $(\theta,t) \in S_{x_1} \times \mathbb{R^{+}}$. For any $\theta \in \Gamma$, let $r(\theta)$ be the radius such that $exp_{x_1}(r (\theta))\in H$. Then it follows from Lemma \ref{Gromov} that $W\subset \{exp_{x_1}(r)|r(\theta) \le r \le 2r(\theta),\  \theta \in \Gamma\}$, and hence
\be\label{vol M}
\vol(W) \le \int_{\Gamma}\int_{r(\theta)}^{2r(\theta)} A(\theta,t)dtd\theta.
\ee
For $r(\theta) \le t \le 2r(\theta) \le 2D_1$, by \eqref{AC}, we have
\[\frac{A(\theta,t)}{t^{m-1}}\le e^{2Lt+Kt^2} \frac{A(\theta,r(\theta))}{r(\theta)^{m-1}},\]
which implies that
\[A(\theta,t) \le e^{4LD_1+4KD_1^2}2^{m-1}A(\theta,r(\theta)).\]
Plugging the above inequality into \eqref{vol M} gives  \[\vol(W)\le e^{4LD_1+4KD_1^2}2^{m-1}\int_{\Gamma} r(\theta)A(\theta,r(\theta))d\theta \le D_12^{m-1} e^{4LD_1+4KD_1^2}  \vol(H^{'}). \] \qed\\

When $W$ is the intersection of $\Omega$ and a ball in $M$, the above lemma implies that

\begin{corollary}
Let $H$ be any hypersurface dividing a convex domain $\Omega$ into two parts $\Omega_1,\Omega_2$. For any ball $B_r(x)$ in $M$, we have
	\begin{equation}
\min(\vol(B_r(x)\cap\Omega_1),\vol(B_r(x)\cap\Omega_2)) \le  2^{m+1}re^{4Ld+4Kd^2}\vol(H\cap(B_{2r}(x))),
	\end{equation}
where $d=diam(\Omega)$, the diameter of $\Omega$. In particular, if $B_r(x)\cap\Omega$ is divided equally by $H$, then
\begin{equation} \label{RV}
\vol(B_r(x)\cap\Omega) \le  2^{m+2}re^{4Ld+4Kd^2}\vol(H\cap B_{2r}(x)) .
	\end{equation}
\proof Put $W_i=B_r(x)\cap\Omega_i$ in the above lemma and use $D_1 \le 2r$ and $H^{'}\subset H\cap B_{2r}(x).$ \qed\\
\end{corollary}

 Now we are ready to prove Theorem \ref{isoperimetric estimate}.\\

\noindent{\it Proof of Theorem \ref{isoperimetric estimate}.} Let $H$ be any hypersurface dividing $M$ into two parts, $\Omega_1$ and $\Omega_2$. We may assume that $\vol(\Omega_1) \le \vol(\Omega_2)$. For any $x\in\Omega_1$, Let $r_x$ be the smallest radius such that
\[\vol(B_{r_x}(x)\cap\Omega_1)=\vol(B_{r_x}(x)\cap\Omega_2)=\frac{1}{2}\vol(B_{r_x}(x)\cap\Omega).\]
By \eqref{RV}, we have,
\begin{equation}\label{CC}
\vol(B_{r_x}(x)\cap\Omega) \le  2^{m+2}r_xe^{4Ld+4Kd^2}\vol(H\cap B_{2r_x}(x)).
\end{equation}
The domain $\Omega_1$ has a covering
\[\Omega_1\subset \cup_{x\in\Omega_1}B_{2r_x}(x).\]
By Vitali Covering Lemma, we can choose a countable family of disjoint balls $B_i=B_{2r_{x_i}}(x_i)$ such that $\cup_iB_{10r_{x_i}}(x_i) \supset \Omega_1.$
So\[\vol(\Omega_1)\le \sum_i \vol(B_{10r_{x_i}}(x_i)\cap\Omega_1).\]
Applying the volume comparison Theorem \ref{volume element comparison} in $\Omega_1$  gives
\[\frac{\vol(B_{10r_{x_i}}(x_i)\cap\Omega_1)}{(10r_{x_i})^m}\le e^{99Kr_{x_i}^2+18Lr_{x_i}}\frac{\vol(B_{r_{x_i}}(x_i)\cap\Omega_1)}{(r_{x_i})^m}.\]
On the other hand, since $\vol(\Omega_1) \le \vol(\Omega_2) $, we have $r_x \le d$ for any $x\in \Omega_1$. Thus,
\begin{align}
\vol(B_{10r_{x_i}} (x_i) \cap\Omega_1)
&\le 10^me^{99Kd^2+18Ld} \vol(B_{r_{x_i}}(x_i)\cap\Omega_1)\nonumber\\ &=2^{-1}10^me^{99Kd^2+18Ld} \vol(B_{r_{x_i}}(x_i)\cap\Omega)\nonumber.
\end{align}
Therefore,
\be\label{vol omega 1}
\vol(\Omega_1) \le 2^{-1}10^me^{99Kd^2+18Ld} \sum_i\vol(B_{r_{x_i}}(x_i)\cap\Omega) .
\ee
Moreover, since the balls $B_i$ are disjoint, \eqref{CC} gives
\be\label{vol H}
\vol(H)\ge \sum_i \vol(B_i \cap H) \ge 2^{-m-2}e^{-4Ld-4Kd^2} \sum_i r_{x_i}^{-1} \vol(B_{r_{x_i}}(x_i)\cap\Omega).
\ee
When $1\le \alpha \le \frac{m}{m-1}$, it follows from \eqref{vol omega 1} and \eqref{vol H} that
\be\label{iso 1}
\begin{split}
\frac{\vol(H)}{\vol(\Omega_1)^{\frac{1}{\alpha}}}   &\ge  \frac{2^{-m-2}e^{-4Ld-4Kd^2}}{(2^{-1}10^me^{99Kd^2+18Ld})^{\frac{1}{\alpha}}}\frac{\sum_i r_{x_i}^{-1}\vol(B_{r_{x_i}}(x_i)\cap\Omega)}{(\sum_i\vol(B_{r_{x_i}}(x_i)\cap\Omega))^{\frac{1}{\alpha}}}  \\
& \ge \frac{2^{-m-2}e^{-4Ld-4Kd^2}}{2^{-1}10^me^{99Kd^2+18Ld}}\frac{\sum_ir_{x_i}^{-1} \vol(B_{r_{x_i}}(x_i)\cap\Omega)}{\sum_i\vol(B_{r_{x_i}}(x_i)\cap\Omega)^{\frac{1}{\alpha}}}  \\
& \ge2^{-2m-1}5^{-m}e^{-22Ld-103Kd^2}\inf_i\frac{r_{x_i}^{-1}\vol(B_{r_{x_i}}(x_i)\cap\Omega)}{\vol(B_{r_{x_i}}(x_i)\cap\Omega)^\frac{1}{\alpha}} \\
&= 2^{-2m-1}5^{-m}e^{-22Ld-103Kd^2}\inf r_{x_i}^{-1}\vol(B_{r_{x_i}}(x_i)\cap\Omega)^{1-\frac{1}{\alpha}}.
\end{split}
\ee
Applying the volume comparison Theorem \ref{volume element comparison} in $\Omega$ gives
\[\frac{\vol(B_d(x_i)\cap\Omega)}{d^m}\le e^{Kd^2+2Ld}\frac{\vol(B_{r_{x_i}}(x_i)\cap\Omega)} {r_{x_i}^m}.\]
Since $1-\frac{1}{\alpha} \ge 0$, and $m(1-\frac{1}{\alpha})-1\le0$,
we derive
\be\label{iso 2}
\al\inf r_{x_i}^{-1} \vol(B_{r_{x_i}}(x_i)\cap\Omega)^{1-\frac{1}{\alpha}} \ge & d^{m(1-\frac{1}{\alpha})-1}\inf r_{x_i}^{-m(1-\frac{1}{\alpha})} \vol(B_{r_{x_i}}(x_i)\cap\Omega)^{1-\frac{1}{\alpha}}\\
\geq & d^{-1}e^{-(Kd^2+2Ld)(1-\frac{1}{\alpha})}\vol(\Omega)^{1-\frac{1}{\alpha}}.
\eal
\ee
From \eqref{iso 1} and \eqref{iso 2}, we conclude that
\[IN_{\alpha}(\Omega)\ge  d^{-1}2^{-2m-1}5^{-m}e^{-(24-\frac{2}{\alpha})Ld-(104-\frac{1}{\alpha})K d^2}\vol(\Omega)^{1-\frac{1}{\alpha}}.\]

On the other hand, when $0<\alpha<1$, similarly to \eqref{iso 1}, we have
\be\label{iso 3}
\begin{split}
\frac{\vol(H)}{\vol(\Omega_1)^{\frac{1}{\alpha}}}&=\frac{\vol(H)}{\vol(\Omega_1) \vol(\Omega_1)^{\frac{1}{\alpha}-1}}    \ge \frac{\vol(H)}{\vol(\Omega_1) \vol(\Omega)^{\frac{1}{\alpha}-1}} \\
	& \ge \frac{2^{-m-2}e^{-4Ld-4Kd^2}}{2^{-1}10^me^{99Kd^2+18Ld}}\frac{\sum_ir_{x_i}^{-1} \vol(B_{r_{x_i}}(x_i)\cap\Omega)}{\sum_i\vol(B_{r_{x_i}}(x_i)\cap\Omega)} \vol(\Omega)^{1-\frac{1}{\alpha}} \\
	& \ge d^{-1}2^{-2m-1}5^{-m}e^{-22Ld-103Kd^2}\vol(\Omega)^{1-\frac{1}{\alpha}} .
	\end{split}
\ee
Taking infimum over $H$ finishes the proof. \qed\\


From Lemma \ref{Cheeger} and Theorem \ref{isoperimetric estimate}, we immediately have the estimate of the first eigenvalue.
\begin{theorem}\label{thm 1}
Let $(M^m,g)$ be an $m$-dimensional closed Riemannian manifold with diameter bounded from above by $D$, and $m\geq 2$. Suppose that \eqref{basic assumption1} and \eqref{basic assumption2} are satisfied. Then
\begin{equation}
	\lambda_1\ge \frac{1}{16}D^{-2}400^{-m}e^{-44LD-206KD^2}:= c_0.
\end{equation}
\end{theorem}

\medskip
To derive the lower bound of higher order eigenvalues, we need to use gradient estimates for eigenfunctions, which in term require a Sobolev inequality. According to section 9 in \cite{Li}, the desired Sobolev inequality follows from the lower bound estimate of $IN_{\frac{m}{m-1}}(M)$.

\begin{definition}[\cite{Li}]
Let	$(M^m,g)$ be an $m$-dimensional compact Riemannian manifold (with or without boundary). For any $\alpha>0$, the Neumann $\alpha$-Sobolev constant of $M$ is defined by
\[SN_\alpha(M)=\inf_{f\in H^{1,1}(M)} \frac{\int_M |\nabla f|}{\{\inf_{k\in \mathbb{R}} \int_M |f-k|^\alpha\}^{\frac{1}{\alpha}}},\]
where $H^{1,1}(M)$ is the Sobolev space.	
	
\end{definition}
As pointed out in \cite{Li}, when $\alpha>\frac{m}{m-1}$, it holds that $IN_{\alpha}(M)=SN_{\alpha}(M)=0$. In general, the relation between $IN_\alpha(M)$ and $SN_\alpha(M)$ is as follows.
\begin{lemma}[section 9 in \cite{Li}]\label{SN and IN}
For any $\alpha >0$,  we have
$$ \min\{{1,2^{1-\frac{1}{\alpha}}}\} IN_\alpha(M) \le SN_\alpha(M)\le \max\{{1,2^{1-\frac{1}{\alpha}}}\}IN_\alpha(M).	$$		
\end{lemma}

Moreover, a lower bound of the Sobolev constant $SN_{\alpha}(M)$ provides a Sobolev inequality. In fact, we have

\begin{lemma}[Corollary 9.9 in \cite{Li}]\label{SN and Sobolev}
Let $(M^m, g)$ be a compact Riemannian manifold (with or without boundary). There exist constants $C_1(\alpha), C_2(\alpha)>0$ depengding only on $\alpha$, such that
\[\int_M |\nabla f|^2 \ge C_1(\alpha)SN_\alpha(M)^2\left(\left(\int_M |f|^{\frac{2\alpha}{2-\alpha}}\right)^\frac{2-\alpha}{\alpha}-C_2(\alpha)\vol(M)^{\frac{(2-2\alpha)}{\alpha}} \int_M |f|^2\right)\]	for all $f\in H^{1,2}(M).$	
\end{lemma}

Then by choosing $\alpha=\frac{m}{m-1}$ for $m\geq 3$ and $\alpha=\frac{4}{3}$ for $m=2$, and combining Lemma \ref{SN and IN}, Lemma \ref{SN and Sobolev}, and Theorem \ref{isoperimetric estimate}, one can get the following Sobolev inequalities.

\begin{corollary} \label{Sobolev inequality}
	Let $(M^m,g)$ be an $m$-dimensional compact Riemannian manifold (with or without boundary). Assume that \eqref{basic assumption1} and \eqref{basic assumption2} are satisfied. Then for any $f\in H^{1,2}(M)$,\\
(1) when $m\ge 3$, we have
	\begin{equation}\label{Sobolev1}
		\int_M |\nabla f|^2 \ge C_1(m) \tilde{C}^2 \vol(M)^{\frac{2}{m}}\left(\left(\int_M |f|^{\frac{2m}{m-2}}\right)^{\frac{m-2}{m}}-C_2(m)\vol(M)^{-\frac{2}{m}}\int_M|f|^2\right),
\end{equation}
where $\tilde{C}=D^{-1}2^{-2m-1}5^{-m}e^{-(22+\frac{2}{m})LD-(103+\frac{1}{m})K D^2}$, and $C_1(m)$ and $C_2(m)$ are dimensional constants;  \\
(2) when $m=2$, one has
\begin{equation}\label{Sobolev2}
		\int_M |\nabla f|^2 \ge \tilde{S_1} \tilde{S}^2 \vol(M)^{\frac{1}{2}}\left(\left(\int_M |f|^{4}\right)^{\frac{1}{2}}-\tilde{S_2}\vol(M)^{-\frac{1}{2}}\int_M|f|^2\right),
\end{equation}
where $\tilde{S}_1$ and $\tilde{S}_2$ are pure constants, and
	$\tilde{S}=D^{-1}2^{-5}5^{-2}e^{-(22+\frac{1}{2})LD-(103+\frac{1}{4})K D^2}.$\\
\end{corollary}	

\begin{remark}\label{constants}
By carefully following the proof of Corollary 9.9 in \cite{Li}, one can check that we may take $C_1(m)=\frac{(m-2)^2}{4(m-1)^2}2^{\frac{2-m}{m(m-1)}}$, $C_2(m)=2^{\frac{2m^3-7m^2+2m+4}{m(m-1)(m-2)}}$, $\tilde{S_1}=3^{-2}2^{-\frac{1}{6}}$, and $\tilde{S_2}=2^{\frac{7}{6}}$.
\end{remark}

\section{Gradient and higher order eigenvalue estimates }	
	
In this section, we use a method in \cite{WZ} (see also \cite{LZZ}) to show the lower bound estimates of high order eigenvalues. Firstly, we prove a gradient estimate of eigenfunctions by Moser iteration.	
	\begin{proposition}\label{prop gradient estimate eigenfunction}
		Let $(M^m,g)$, $m\ge3$, be an $m$-dimensional closed Riemannian manifold. Suppose that \eqref{basic assumption1} and \eqref{basic assumption2} are satisfied. Let $\lambda$ be an eigenvalue of the Laplace operator, and $u$ an eigenfunction  satisfing $\Delta u=-\lambda u$. Then we have the following gradient estimate.
		\begin{equation}
			|\nabla u|^2\le 2^m\left(\frac{m}{m-2}\right)^\frac{m(m-2)}{2}\left(\frac{3\lambda+2K+2L^2+C_2}{C_1}\right)^\frac{m}{2}(\lambda+L^2)\vol(M)^{-1}\int_M u^2,
		\end{equation}
	where $C_1=C_1(m)\tilde{C}^2$, and $C_2=C_1(m)\tilde{C}^2C_2(m)$ with $C_1(m),\ C_2(m)$, $\tilde{C}$ the constants in \eqref{Sobolev1}.

In particular, when $||u||_{L^2}=1$, we have
			\begin{equation}\label{gradient estimate eigenfunction}
		|\nabla u|^2\le 2^m\left(\frac{m}{m-2}\right)^\frac{m(m-2)}{2}\left(\frac{3\lambda+2K+2L^2+C_2}{C_1}\right)^\frac{m}{2}(\lambda+L^2)\vol(M)^{-1}.
			\end{equation}
\end{proposition}	
	
\proof Let $v=|\nabla u|^2+L^2u^2$. The Bochner formula and assumptions \eqref{basic assumption1} and \eqref{basic assumption2} induce that
\begin{align}
	\Delta v &=2|Hess\,u|^2+2<\nabla\Delta u,\nabla u>+2Ric(\nabla u,\nabla u)+2L^2u\Delta u+2L^2|\nabla u|^2 \nonumber\\
		&\ge 2|Hess\,u|^2-2\lambda|\nabla u|^2-2K|\nabla u|^2-2f_{ij}u_iu_j-2L^2\lambda u^2+2L^2|\nabla u|^2\nonumber\\
		&=2|Hess\,u|^2-2\lambda v +(2L^2-2K)|\nabla u|^2-2f_{ij}u_iu_j\nonumber\\
		&\ge 2u_{ij}^2-2(\lambda+K) v-2f_{ij}u_iu_j.\nonumber	
\end{align}
Multiple both sides above by $v^{p-1}$, $p\geq 2$, and take integrals over $M$. Notice that
\begin{align}
	\int_Mv^{p-1}\Delta v&=-\int_M<\nabla v^{p-1},\nabla v>=-\int_M(p-1)v^{p-2}<\nabla v,\nabla v>\nonumber\\
	&=-(p-1)\int_Mv^{p-2}|\nabla v|^2=-
	\frac{4(p-1)}{p^2}\int_M|\nabla v^{\frac{p}{2}}|^2\nonumber.
\end{align}
Hence, we have
\begin{equation}\label{multiply v^p-1}
\frac{4(p-1)}{p^2}\int_M|\nabla v^{\frac{p}{2}}|^2\le -2\int_M u_{ij}^2v^{p-1}+2(\lambda+K) \int_Mv^p+2\int_Mf_{ij}
	u_iu_jv^{p-1}.
\end{equation}
For the third term on the right hand side above, integrating by part yields
\begin{align}\label{int by parts}
	2\int_Mf_{ij}u_iu_jv^{p-1}&=-2\int_Mf_i(u_iu_jv^{p-1})_j\nonumber\\
	&=\underbrace{-2\int_Mf_iu_{ij}u_jv^{p-1}}_{I}\underbrace{-2\int_Mf_iu_iu_{jj}v^{p-1}}_{II}\underbrace{-2\int_Mf_iu_iu_j(p-1)v^{p-2}v_j}_{III}.
\end{align}
For $I$ above, using Cauchy-Schwarz inequality and the bound of $|\nabla f|$ gives
	\begin{align}
	I=-2\int_Mf_iu_{ij}u_jv^{p-1}& \le 2\int_M(u_{ij}^2v^{p-1}+\frac{1}{4}f_i^2u_j^2v^{p-1})\nonumber\\
	&=2\int_Mu_{ij}^2v^{p-1}+\frac{1}{2}\int_Mf_i^2u_j^2v^{p-1}\nonumber\\
	&\le 2\int_Mu_{ij}^2v^{p-1}+\frac{L^2}{2}\int_Mv^{p}.\nonumber
\end{align}
Next, noticing that $v\ge2L|\nabla u||u|$, we have
	$$II=-2\int_Mf_iu_iu_{jj}v^{p-1}=2\lambda\int_Mf_iu_iuv^{p-1}\le2\lambda L \int_M|\nabla u||u|v^{p-1}\le\lambda \int_M v^p.$$
Finally, by applying Cauchy-Schwarz inequality inequality again to $III$, we deduce
\begin{align}
	III=-2\int_Mf_iu_iu_j(p-1)v^{p-2}v_j&\le2(p-1)L\int_M|\nabla u|^2|\nabla v|v^{p-2}\le2(p-1)L\int_M|\nabla v|v^{p-1}\nonumber\\
	&\le 2(p-1)L(\frac{1}{4\epsilon_1}\int_Mv^p+\epsilon_1\int_M|\nabla v|^2v^{p-2})\nonumber\\
	&=\frac{(p-1)L}{2\epsilon_1}\int_Mv^p+\frac{8(p-1)L\epsilon_1}{p^2}\int_M
	|\nabla v^{\frac{p}{2}}|^2,\nonumber
\end{align}
where $\epsilon_1> 0$ is any constant. Thus, by combining the above estimates in \eqref{multiply v^p-1}, we arrive at
	$$(\frac{4(p-1)}{p^2}-\frac{8(p-1)L\epsilon_1}{p^2})\int_M|\nabla v^{\frac{p}{2}}|^2\le(3\lambda+\frac{L^2}{2}+ 2K+\frac{(p-1)L}{2\epsilon_1})\int_Mv^p.$$
Assume for now that $L>0$. Then, by choosing $\epsilon_1=\frac{1}{4L}$ and noticing that $\frac{2(p-1)}{p^2}\ge\frac{1}{p}$ for $p\geq 2$, one gets
\be\label{L>0}
\int_M|\nabla v^{\frac{p}{2}}|^2\le p^2(3\lambda+2L^2+ 2K)\int_Mv^p.
\ee
If $L=0$, then $f$ is a constant, and from \eqref{multiply v^p-1} we conclude that
	\[\int_M|\nabla v^{\frac{p}{2}}|^2\le \frac{1}{2}p^2(\lambda+K)\int_Mv^p,\]
which is better than \eqref{L>0}. Therefore, we always have
\be\label{L}
\int_M|\nabla v^{\frac{p}{2}}|^2\le p^2(3\lambda+2L^2+ 2K)\int_Mv^p.
\ee
Recall the Sobolev inequality \eqref{Sobolev1},
	\begin{align} \label{22222}
		\int_M |\nabla f|^2 \ge C_1\vol(M)^{\frac{2}{m}}\left(\int_M|f|^{\frac{2m}{m-2}}\right)^{\frac{m-2}{m}}-C_2\int_M|f|^2
	\end{align}	
for all $f\in H^{1,2}(M)$, where $C_1=C_1(m)\tilde{C}^2$, and  $C_2=C_1(m)\tilde{C}^2C_2(m).$
Putting $f=v^{\frac{p}{2}}$ and using \eqref{L} yield
	$$\left(\int_Mv^{\frac{pm}{m-2}}\right)^{\frac{m-2}{m}}\le p^2\left(\frac{3\lambda+2K+2L^2+C_2}{C_1\vol(M)^{\frac{2}{m}}}\right)\int_Mv^p.$$
Denote $Q=\frac {3\lambda+2K+2L^2+C_2}{C_1\vol(M)^{\frac{2}{m}}}$ \ for convenience. The inequality above means that
	$$||v||_{\frac{pm}{m-2}}\le(p^2Q)^\frac{1}{p}||v||_p$$ for all $p\ge2$.

Setting $\beta=\frac{m}{m-2},\ p=2\beta^j $ for $ j=0,\ 1,\ 2,\ ...\, ,$ it implies that
	$$||v||_{2\beta^{j+1}}\le2^{\frac{1}{\beta^j}}\beta^{\frac{j}{\beta^j}}Q^{\frac{1}{2\beta^j}}||v||_{2\beta^j}.$$
Iterating this estimate, we conclude that
	$$||v||_{2\beta^{j+1}}\le2^{\sum_{l=0}^{j}\frac{1}{\beta^l}}\beta^{\sum_{l=0}^{j}\frac{l}{\beta^l}}Q^{\sum_{l=0}^{j}\frac{1}{2\beta^l}}||v||_2.$$	Letting $j\to\infty$, we obtain
	$$||v||_{\infty}\le2^{\frac{m}{2}}\left(\frac{m}{m-2}\right)^\frac{m(m-2)}{4}Q^{\frac{m}{4}}||v||_2.$$
Notice that $\int_M v^2\le||v||_{\infty}\int_M v$. Therefore, the above estimate reduces to
	$$\max\limits_{M} v\le2^m\left(\frac{m}{m-2}\right)^\frac{m(m-2)}{2}Q^{\frac{m}{2}}\int_Mv.$$
This finishes the proof, since $$\int_M v=\int_M(|\nabla u|^2+L^2u^2)=(\lambda+L^2)\int_Mu^2.$$ \qed\\

When $m=2$, by using the Sobolev inequality \eqref{Sobolev2} instead of \eqref{Sobolev1}, one can similarly obtain the following gradient estimate for $u$.
\begin{proposition}\label{prop gradient estimate eigenfunction m=2}
If $(M, g)$ is a Riemann surface, $u$ is an eigenfunction associated to eigenvalue $\lambda$, and \eqref{basic assumption1} and \eqref{basic assumption2} are satisfied, then
\[|\nabla u|^2 \le 2^8\left(\frac{3\lambda+2K+2L^2+S_2}{S_1}\right)^2(\lambda +L^2)\vol(M)^{-1}\int_M u^2,\]
where $S_1=\tilde{S_1}\tilde{S}^2$, and $S_2=\tilde{S_1}\tilde{S}^2\tilde{S_2}$ with $\tilde{S_1},\ \tilde{S_2}$, $\tilde{S}$ the constants in \eqref{Sobolev2}.
\end{proposition}

Next, we prove a similar gradient estimate for linear combinations of eigenfunctions.

\begin{proposition}\label{prop gradient estimate combination}
Let $(M^m, g_{ij})$ be an $m$-dimensional closed Riemannian manifold satisfying \eqref{basic assumption1} and \eqref{basic assumption2}. Let $\phi_j$ be a normalized eigenfunction associated to $\lambda_j$, $j=1,\,2,\,\ ...,\ k$ i.e., $\Delta \phi_j=-\lambda_j \phi_j$ and $\int_M |\phi_j|^2 dV=1$. Then for any sequence of real numbers $b_j,\ j=1,\ 2,\ ...,\ k,$ with $\sum_{j=1}^{k}b_j^2 \le 1$, the linear combination $w=\sum_{i=1}^{k}b_j\phi_j$ satisfies that, for $m\geq 3$,
\begin{equation}\label{gradient estimate combination}	
		|\nabla w|^2 +L^2w^2 \le 2^m \left(\frac{m}{m-2}\right)^\frac{m(m-2)}{2}\left(\frac{6\lambda_k+2K+2L^2+C_2}{C_1}\right)^{\frac{m}{2}}(\lambda_k+L^2)\vol(M)^{-1},
	\end{equation}
and for $m=2$,
\be \label{gradient estimate combination m=2}
|\nabla w|^2+L^2w^2 \le 2^8\left(\frac{6\lambda_k+2K+2L^2+S_2}{S_1}\right)^2(\lambda_k+L^2)\vol(M)^{-1},
\ee
where $C_1,\ C_2,\  S_1,\  S_2$ are constants in Propositions \ref{prop gradient estimate eigenfunction} and \ref{prop gradient estimate eigenfunction m=2}.

\end{proposition}

\proof Here, we only present the proof of \eqref{gradient estimate combination}. The proof of \eqref{gradient estimate combination m=2} is similar by using \eqref{Sobolev2} instead of \eqref{Sobolev1}. First of all, since $\lambda_k>0$, we can write
		 $$\Delta w=-\sum_{j=1}^{k}\lambda_jb_j\phi_j=-\lambda_k \eta,$$ where $\displaystyle \eta=\sum_{j=1}^{k}\frac{\lambda_j}{\lambda_k}b_j\phi_j.$\\

Let $v=|\nabla w|^2+L^2w^2$. Then
		  	\begin{align}	
		\Delta v &=2|Hess\,w|^2+2<\nabla \Delta w,\nabla w>+2Ric(\nabla w,\nabla w)+2L^2w\Delta w +2L^2|\nabla w|^2\nonumber\\	
		&\ge 2w_{ij}^2-2\lambda_k\eta_iw_i-2K|\nabla w|^2-2f_{ij}w_iw_j-2L^2\lambda_k\eta w.\nonumber\\
		&\ge 2w_{ij}^2-2\lambda_k\eta_iw_i-2Kv-2f_{ij}w_iw_j-2L^2\lambda_k\eta w.\nonumber
	\end{align}			
		Multiplying both sides by $v^{p-1}, \  p \ge 2$, and integrating over $M$ give
			\begin{equation} \label{222}
				\begin{split}
					\frac{4(p-1)}{p^2}\int_M|\nabla v^{\frac{p}{2}}|^2 &\le -2\int_M w_{ij}^2v^{p-1}+2\lambda_k\int_M \eta_iw_iv^{p-1}\\
					&+2K\int_Mv^p+2\int_M f_{ij}w_iw_jv^{p-1} +2\lambda_k L^2\int_M \eta wv^{p-1}.
				\end{split}
	\end{equation}
Using H\"older inequality yields
\begin{equation}\label{combination1}
		2\lambda_k\int_M \eta_iw_iv^{p-1} \le 2\lambda_k  \int_M |\nabla \eta|v^{p-\frac{1}{2}} \le 2\lambda_k \left( \int_M v^p \right)^{\frac{p-\frac{1}{2}}{p}} \left( \int_M |\nabla\eta|^{2p}\right)^{\frac{1}{2p}}.
\end{equation}
Notice that the coefficients in $\nabla \eta$ satisfy $\sum_{j=1}^{k}(\frac{\lambda_j}{\lambda_k}b_j)^2 \le \sum_{j=1}^{k}b_j^2 \le 1$ and $\int_M v^p \ge \int_M |\nabla w|^{2p}$. Thus,
\be\label{combination2}
\int_M |\nabla \eta|^{2p} \le \max \limits_{b_1,\dots,b_k} \int_M v^p.
\ee
By combining \eqref{combination1} and \eqref{combination2}, we obtain
\be\label{max1}
2\lambda_k\int_M \eta_iw_iv^{p-1} \le 2\lambda_k \max \limits_{b_1,\dots,b_k} \int_M v^p.
\ee
Here and in the rest of the proof, the maximum is taken for all real numbers $b_1,\cdots,b_k$ such that $\sum_{j=1}^k b_j^2\leq1$.

Similarly, for the last term of \eqref{222}, we have		
\be
\al
2\lambda_kL^2 \int_M \eta wv^{p-1} \le& 2 \lambda_kL\int_M |\eta|v^{p-\frac{1}{2}} \le 2\lambda_kL \left( \int_M v^p \right)^{\frac{p-\frac{1}{2}}{p}} \left( \int_M |\eta|^{2p}\right)^{\frac{1}{2p}}\\
\leq &2\lambda_k \max \limits_{b_1,\dots,b_k} \int_M v^p
\eal
\ee

Finally, we need to deal with the fourth  term on the right hand side of \eqref{222}. Using integration by parts gives
\begin{equation}\label{22}
2\int_M f_{ij}w_iw_jv^{p-1}=\underbrace{-2\int_M f_{i}w_{ij}w_jv^{p-1}}_{I}\underbrace{-2\int_M f_{i}w_iw_{jj}v^{p-1}}_{II}\underbrace{-2\int_M f_{i}w_iw_j(p-1)v^{p-2}v_j}_{III}.
\end{equation}
Using Cauchy-Schwarz inequality and the bound of $|\d f|$, we have
\[I=-2\int_M f_{i}w_{ij}w_jv^{p-1}\le 2\int_M w^2_{ij}v^{p-1}+\frac{L^2}{2}\int_M v^p \le 2\int_M w^2_{ij}v^{p-1} +\frac{L^2}{2} \max \limits_{b_1,\dots,b_k} \int_M v^p,\]
\[II=-2\int_M f_{i}w_iw_{jj}v^{p-1}\le 2\lambda_k \int_M |\nabla f| |\nabla w| |\eta| v^{p-1} \le 2\lambda_k L \int_M |\eta|v^{p-\frac{1}{2}} \le 2\lambda_k\max \limits_{b_1,\dots,b_k} \int_M v^p,\]
and
\begin{align}
		III=-2\int_M f_{i}w_iw_j(p-1)v^{p-2}v_j &\le 2(p-1)L \int_M |\nabla w |^2v^{p-2}|\nabla v| \le 2(p-1)L \int_M v^{p-1}|\nabla v|\nonumber\\
		& \le 2(p-1)L \left(\frac{1}{4\varepsilon_2}\int_M v^p +\varepsilon_2 \int_M v^{p-2}|\nabla v|^2 \right)\nonumber\\
		&= \frac{(p-1)L}{2\varepsilon_2}\max_{b_1,\dots,b_k}\int_M v^p +\frac{8(p-1)L\varepsilon_2}{p^2} \int_M |\nabla v^{\frac{p}{2}}|^2, \nonumber
			\end{align}	
where $\varepsilon_2 >0$ is arbitrary constant. Hence, plugging the estimates above in \eqref{222} asserts that	
	\[	\left(\frac{4(p-1)}{p^2}-\frac{8(p-1)L\varepsilon_2}{p^2} \right ) \int_M |\nabla v^{\frac{p}{2}}|^2 \le \left(6\lambda_k +\frac{L^2}{2}+2K+\frac{(p-1)L}{2\varepsilon_2}\right) \max \limits_{b_1,\dots,b_k} \int_M v^p .\]
Choosing $\varepsilon_2=\frac{1}{4L}$, it follows that
\be\label{max2}
\max_{b_1,\dots,b_k}\int_M|\nabla v^{\frac{p}{2}}|^2 \le p^2 \left(6\lambda_k+2K+2L^2\right) \max \limits_{b_1,\dots,b_k} \int_M v^p.
\ee
Again, by \eqref{max2} and the Sobolev inequality \eqref{Sobolev1}, we have
			\begin{equation}	
		\max \limits_{b_1,\dots,b_k}  \left(\int_M v^{\frac{pm}{m-2}} \right)^{\frac{m-2}{m}} \le p^2
		\left(\frac{6\lambda_k+2K+2L^2+C_2}{C_1\vol(M)^{\frac{2}{m}}}\right) \max \limits_{b_1,\dots,b_k} \left(\int_M v^p \right).
			\end{equation}
Denoting $Q=\frac{6\lambda_k+2K+2L^2+C_2}{C_1\vol(M)^{\frac{2}{m}}}$ and using Moser iteration as in Proposition \ref{prop gradient estimate eigenfunction}, it follows that
		\[\max \limits_{b_1,\dots,b_k} ||v||_{\infty} \le 2^{\frac{m}{2}}  \left(\frac{m}{m-2} \right)^{\frac{m(m-2)}{4}}Q^{\frac{m}{4}} \max \limits_{b_1,\dots,b_k} ||v||_2 .\]
Square both sides above and notice that
		\[ \max \limits_{b_1,\dots,b_k} \int_M v^2 \le \max \limits_{b_1,\dots,b_k} ||v||_{\infty} \max \limits_{b_1,\dots,b_k}\int_M v. \]
Thus, we get \\
			\begin{equation} \label{3}	
		\max \limits_{b_1,\dots,b_k} ||v||_{\infty} \le 2^m  \left(\frac{m}{m-2} \right)^{\frac{m(m-2)}{2}}Q^{\frac{m}{2}} \max \limits_{b_1,\dots,b_k}\int_M v .
	\end{equation}

On the other hand, since $\phi_1,\ \phi_2,\ \dots,\ \phi_k$ are orthonormal, we have
		\[
		\begin{split}
		\int_M v &=\int_M (|\nabla w|^2+L^2w^2)=-\int_M w\Delta w+L^2\int_M w^2 \\ &=\int_M(\sum_{j=1}^{k}b_j\phi_j)(\sum_{i=1}^{k} \lambda_i b_i \phi_i)+L^2\int_M(\sum_{j=1}^{k}b_j\phi_j)^2 \\
			&=\sum_{j=1}^{k} \lambda_j b_j^2+L^2\sum_{j=1}^{k}b_j^2 \le (\lambda_k+L^2)\sum_{j=1}^{k}b_j^2 \le \lambda_k+L^2.
		\end{split}
		\]
This, together with \eqref{3}, completes the  proof.\qed\\
		
The above gradient estimate for linear combinations of eigenfunctions allows us to derive the arithmetic inequality of the eigenvalues below.	
	
\begin{lemma}\label{lem combination eigenvalue}
Under the same assumptions and notations as in Proposition \ref{prop gradient estimate combination}, we have for $m\geq 3$,
\begin{equation}\label{combination eigenvalue}
\lambda_1+\lambda_2+...+\lambda_k\le \\
m 2^m \left(\frac{m}{m-2}\right)^\frac{m(m-2)}{2}\left(\frac{6\lambda_k+2K+2L^2+C_2}{C_1}\right)^{\frac{m}{2}}(\lambda_k+L^2),
\end{equation}
and for $m=2$,
\be\label{combination eigenvalue m=2}
\lambda_1+\lambda_2+...+\lambda_k \le 2^9\left(\frac{6\lambda_k+2K+2L^2+S_2}{S_1}\right)^2(\lambda_k+L^2).
\ee
\end{lemma}

\proof We only prove \eqref{combination eigenvalue} by using \eqref{gradient estimate combination}. The proof of \eqref{combination eigenvalue m=2} follows similarly from \eqref{gradient estimate combination m=2}.

If $k\le m$, the conclusion follows immediately from Proposition \ref{prop gradient estimate eigenfunction} by integrating both sides of \eqref{gradient estimate eigenfunction} for each $\phi_j$, $j=1,2,\cdots,k$.

When $k> m$, for each $x\in M$, we can find an orthogonal matrix $(a_{ij})_{k\times k}$ such that
$$\varphi_i=\sum_{j=1}^{k}a_{ij}\phi_j,i=1,\ 2,\ \dots,\ k$$ satisfy that
$$\nabla_l\varphi_i(x)=0,\ l=1,\ 2,\ \dots,\ m,\ m+1\le i \le k. $$
Indeed, since the rank of the matrix
	\begin{equation}
	J=\begin{pmatrix}
		\nabla_1\phi_1&\dots&\nabla_1\phi_k\\
		\vdots& &	\vdots\\
		\nabla_m\phi_1&\dots&\nabla_m\phi_k
	\end{pmatrix}
\end{equation}
is no more than $m$, there are $k-m$ linearly independent solutions of $J\vec{x}=\vec{0}$, and then Schmidt orthogonalization gives $(a_{ij})$.

Thus, we derive from Proposition \ref{prop gradient estimate combination} that
$$|\nabla \phi_1|^2+...+|\nabla \phi_k|^2=|\nabla \varphi_1|^2+...+|\nabla \varphi_k|^2=|\nabla \varphi_1|^2+...+|\nabla \varphi_m|^2 $$
$$\le m 2^m \left(\frac{m}{m-2}\right)^\frac{m(m-2)}{2}\left(\frac{6\lambda_k+2K+2L^2+C_2}{C_1}\right)^{\frac{m}{2}}(\lambda_k+L^2)\vol(M)^{-1}.$$ \qed\\
Thus, integrating both sides gives Lemma \ref{lem combination eigenvalue}.

\begin{remark}
Notice that the above Lemma cannot be deduced directly from Propositions \ref{prop gradient estimate eigenfunction} and \ref{prop gradient estimate eigenfunction m=2}, which will enlarge the coefficient $m$ on the right hand side of \eqref{combination eigenvalue} and \eqref{combination eigenvalue m=2} to be $k$.
\end{remark}

From \eqref{combination eigenvalue} and \eqref{combination eigenvalue m=2}, in order to get a lower bound of $\lambda_k$, we only need the following lemma.

\begin{lemma}[\cite{WZ}]\label{lem WZ}
	For $0\le\lambda_1\le\lambda_2\le...\le\lambda_k\le...$, if the inequality
	\begin{equation}	
		\lambda_1+\lambda_2+...+\lambda_k\le C_3\lambda_k^{\frac{m}{2}+1}
	\end{equation}
holds for any $k\ge 1$, then ones has
	\begin{equation}
		\lambda_k\ge C_4k^{\frac{2}{m}},
	\end{equation}
	\\
	where
	$$C_4=min\left\{\lambda_1,\ \left(\frac{m}{C_3(m+2)}\right)^{\frac{2}{m}}\right\},$$
	and $m\ge 1$ is an integer.
\end{lemma}

Now we can see that a lower bound of $\lambda_k$ follows immediately from Theorem \ref{thm 1}, Lemma \ref{lem combination eigenvalue} and Lemma \ref{lem WZ}.

\begin{theorem}\label{thm n}
Assume that $(M^m,g)$ is an $m$-dimensional closed Riemannian manifold such that \eqref{basic assumption1} and \eqref{basic assumption2} are satisfied. Let $c_0$ be the lower bound of $\lambda_1$ in Theorem \ref{thm 1}. Then\\
(1) for $m\geq 3$,
\begin{equation}\label{lambda k lower bound}
			\lambda_k\ge c_1k^{\frac{2}{m}},\ \forall k\geq 2,
		\end{equation}
		where $c_1=min\left\{c_0,\ \left(\frac{m}{C_5(m+2)}\right)^{\frac{2}{m}}\right\},$ and \\ $C_5=m2^m\left(\frac{m}{m-2}\right)^\frac{m(m-2)}{2}c_0^{-(\frac{m}{2}+1)}\left(\frac{6c_0+2K+2L^2+C_2}{C_1}\right)^\frac{m}{2}(c_0+L^2);$\\
(2) for $m=2$,
\be\label{lambda k lower bound m=2}
\lambda_k \ge c_2k^{\frac{1}{2}},\ \forall k\geq 2,
\ee
where $c_2=\min\left\{c_0,\ \left(\frac{2}{3C_6}\right)^{\frac{1}{2}}\right\}$, and $C_6=2^9c_0^{-3}\left(\frac{6c_0+2K+2L^2+S_2}{S_1}\right)^2\left(c_0+L^2\right). $

\end{theorem}
\proof To prove \eqref{lambda k lower bound}, from Lemma \ref{lem combination eigenvalue}, we have
$$\lambda_1+\lambda_2+...+\lambda_k\le \lambda_k^{\frac{m}{2}+1}m 2^m \left(\frac{m}{m-2}\right)^\frac{m(m-2)}{2}\left(\frac{6+\frac{2K+2L^2+C_2} {\lambda_k}}{C_1}\right)^{\frac{m}{2}}(1+\frac{L^2}{\lambda_k}).$$
Since $\lambda_k\geq \lambda_1\geq c_0$, it follows that
\begin{equation}
	\lambda_1+\lambda_2+...+\lambda_k\le C_5\lambda_k^{\frac{m}{2}+1}.
\end{equation}
From Lemma \ref{lem WZ}, we can easily  get the conclusion.

The proof of \eqref{lambda k lower bound m=2} is similar. \qed

\begin{remark}
Recall that the constants $C_1$, $C_2$, $S_1$, and $S_2$ have explicit expressions according to Corollary \ref{Sobolev inequality} and Remark \ref{constants}. Thus, the lower bound of $\lambda_k$ in the above theorem can also be expressed explicitly.\\
\end{remark}

\section*{Acknowledgements}

Research is partially supported by NSFC Grant No. 11971168, Shanghai Science and Technology Innovation Program Basic Research Project STCSM 20JC1412900, and Science and Technology Commission of Shanghai Municipality (STCSM) No. 18dz2271000.

	\end{document}